\documentclass[reqno,10pt]{article}
\usepackage{amsmath}
\usepackage{graphicx,color}

\newcommand{\be}{\begin{equation}}
\newcommand{\ee}{\end{equation}}
\newenvironment{pf}{\noindent{\it
Proof}.\enspace}{\rule{2mm}{2mm}\medskip}

\newcommand{\R}{\mathbb{R}}
 \newcommand{\Rn}{\mathbb{R}^n}

\newcommand{\E}{\mathbb{E}}

\newcommand{\dyle}{\displaystyle}

\renewcommand{\a }{\alpha }
\renewcommand{\b }{\beta }

\newcommand{\D }{\Delta }

\renewcommand{\l }{\lambda }
\renewcommand{\L }{\Lambda }

\newcommand{\n }{\nabla }
\newcommand{\var }{\varphi }
\renewcommand{\o }{\omega }

\newcommand{\cN}{{\cal N}}

\newcommand{\Wan}{W^{s,2}(\Rn)}

\newcommand{\intn}{\int_{\Rn}}

\newcommand{\intN}{\int_{\R^n}}

\newcommand{\wt}{\widetilde}

\newcommand{\bu}{{\bf u}}
\newcommand{\bv}{{\bf v}}

\newcommand{\bw}{{\bf w}}
\newcommand{\bo}{{\bf 0}}

\newcommand{\bh}{{\bf h}}

\def\diver{\mathop{\rm div}\nolimits}
\def\tr{\mathop{\rm Tr}\nolimits}
\def\ext{\rm Ext_s}
\def\div{\mathop{\rm div}\nolimits}

\numberwithin{equation}{section}

\newtheorem{Theorem}{Theorem}

\newtheorem{Lemma}[Theorem]{Lemma}
\newtheorem{Proposition}[Theorem]{Proposition}
\newtheorem{Definition}[Theorem]{Definition}
\newtheorem{remark}[Theorem]{Remark}
\newtheorem{remarks}[Theorem]{Remarks}
\newtheorem{example}[Theorem]{Example}
\newtheorem{examples}[Theorem]{Examples}
\newenvironment{Remark}{\begin{remark} \rm}{\rule{2mm}{2mm}\end{remark}}

\numberwithin{Theorem}{section}

        \usepackage{latexsym}
        \usepackage{amsmath}
        \usepackage{amsfonts}
        \usepackage{amssymb}

\begin{document}
\title{
\vspace{0.5in} {\bf\Large  Ground states of some coupled nonlocal fractional dispersive PDEs}}
\author{{\bf\large Eduardo Colorado}\footnote{Partially supported
by the Ministry of Economy and Competitiveness of Spain and FEDER under
Research Project MTM2016-80618-P, and by the INdAM - GNAMPA Project 2017
``Teoria e modelli non locali".}\hspace{2mm}
{\bf\large}\vspace{1mm}\\
{\it\small Departamento de Matem\'aticas, Universidad Carlos III de Madrid}\\
{\it\small Avda. Universidad 30, 28911 Legan\'es (Madrid), Spain.}\\
{\it \small \& Instituto de Ciencias Matem\'aticas, ICMAT(CSIC-UAM-UC3M-UCM)}\\
{\it \small C/Nicol\'as Cabrera 15, 28049 Madrid,
Spain}\\
{\it\small e-mail: eduardo.colorado@uc3m.es, \&
eduardo.colorado@icmat.es}}
\date{}

\maketitle

{\it {\small  Dedicated to the memory of Anna Aloe}}

\

\begin{center}
{\bf\small Abstract}

\vspace{3mm} \hspace{.05in}\parbox{4.5in} {\small We show the
existence of ground state solutions to the following stationary system coming from some coupled fractional dispersive equations such as:
nonlinear fractional  Schr\"odinger (NLFS) equations (for dimension $n=1,\, 2,\, 3$) or NLFS and fractional Korteweg-de Vries equations  (for $n=1$),
$$
\left \{
\begin{array}{ll}
(-\D)^{s} u+ \l_1 u &=  u^{3}+\b uv,\quad u\in \Wan\\
(-\D)^{s}  v + \l_2 v &=  \frac 12 v^{2}+\frac 12 \b u^2,\quad v\in \Wan,
\end{array} \right.
$$
where $\l_j>0$, $j=1,2$, $\b\in
\mathbb{R}$, $n=1,\, 2,\, 3$, and $\frac n4< s<1$.
Precisely, we  prove the existence of a positive radially symmetric ground state for any $\b>0$.}
\end{center}

\noindent
{\it \footnotesize 2010 Mathematics Subject Classification}. {\scriptsize 49J40, 35Q55, 35Q53, 35B38, 35J50}.\\
{\it \footnotesize Key words}. {\scriptsize  Nonlinear Fractional Schr\"odinger Equations, Fractional Korteweg-de Vries equations,
Variational Methods, Critical Point Theory, Ground States.}

\maketitle

\section{Introduction}
In this paper we  study the existence of ground state solutions to the following stationary system coming from some coupled nonlocal fractional dispersive
equations such as:
nonlinear fractional  Schr\"odinger (NLFS) equations (for dimension $n=1,\, 2,\, 3$) or NLFS and fractional Korteweg-de Vries equations (FKdV) (for $n=1$)
\begin{equation}\label{eq:Main}
\left \{
\begin{array}{ll}
(-\D)^{s} u+ \l_1 u &=  u^{3}+\b uv,\quad u\in \Wan\\
(-\D)^{s}  v + \l_2 v &=  \frac 12 v^{2}+\frac 12 \b u^2,\quad v\in \Wan,
\end{array} \right.
\end{equation}
where $\Wan$  denotes the fractional Sobolev space,  $n=1,\, 2,\, 3$.
$\l_j>0$, $j=1,2$, the coupling factor $\b\in \mathbb{R}$,
and the fraction $\frac n4< s<1$.

The associated critical Sobolev exponent is defined by $2^*_s= \dfrac{2n}{n-2s}$ if
 $n>2s$, and $2^{*}_s=\infty $ if $n\le 2s$. As a consequence, since $\frac n4<s<1$ we have that $2^*_s>4$.

It is well known that the fractional Laplacian $(-\D)^s$, $0<s<1$, is a
nonlocal diffusive type operator. It arises in several
 physical phenomena like flames propagation and chemical reactions in liquids, population dynamics, geophysical fluid dynamics,
 in probability, American options in finance, in $\alpha$-stable L\'evy processes, etc; see for instance
 \cite{Applebaum,Bertoin,Cont-Tankov}.

In the one-dimensional case, when $s=1$,  \eqref{eq:Main} comes from the following system of coupled nonlinear Sch\"odinger (NLS) and
Korteweg-de Vries (KdV)
equations
\begin{equation}\label{NLS-KdV}
\left\{\begin{array}{rcl}
if_t + f_{xx} + |f|^2f+ \b fg & = &0\\
g_t +g_{xxx} +gg_x  + \frac 12\b(|f|^2)_x  & = & 0,
\end{array}\right.
\end{equation}
where $f=f(x,t)\in \mathbb{C}$ while $g=g(x,t)\in \mathbb{R}$, and
$\b\in \mathbb{R}$ is the real coupling coefficient. System \eqref{NLS-KdV}
appears in phenomena of interactions between short and long
dispersive waves, arising in fluid mechanics, such as  the
interactions of capillary - gravity water waves. Indeed, $f$
represents the short-wave, while $g$ stands for the long-wave. For more details, see for instance
\cite{aa,cl,fo} and the references therein.

Looking for ``traveling-wave" solutions, namely
solutions to \eqref{NLS-KdV} of the form
$$
(f(x,t),g(x,t))=\left(e^{i\o t}
e^{i\frac c2 x}u(x-ct),v(x-ct)\right)\quad\mbox{with}\quad u,\,
v\quad\mbox{real functions,}
$$
 and choosing  $\l_1=\o+\frac{c^2}{4}$,
$\l_2=c$, one finds that $u,\, v$ solve the following problem
\begin{equation}\label{NLS-KdV2}
\left\{\begin{array}{rcl}
-u'' +\l_1 u & = & u^3+\beta uv \\
-v'' +\l_2 v & = & \frac 12 v^2+\frac 12\beta u^2.
\end{array}\right.
\end{equation}
This system \eqref{NLS-KdV2} has been studied, among others, in \cite{aa,ab,c2,c3,dfo1,dfo2}.
Also, note that system \eqref{NLS-KdV2} corresponds to system \eqref{eq:Main} when $s=1$ and $n=1$.

On the other hand, for $n=2,\, 3$, and $s=1$, system \eqref{eq:Main} corresponds to
\eqref{NLS-KdV2}
\begin{equation}\label{NLS-KdV2-n}
\left\{\begin{array}{rcl}
-\D u +\l_1 u & = & u^3+\beta uv \\
-\D v +\l_2 v & = & \frac 12 v^2+\frac 12\beta u^2,
\end{array}\right.
\end{equation}
for which the existence of bound and ground states have been studied in \cite{c2,c3}.
We observe that system \eqref{NLS-KdV2-n} can be seen as a stationary version of a time dependent coupled NLS system
 when one looks for solitary wave solutions, and  $(u,v)$ are the corresponding standing waves solutions of \eqref{NLS-KdV2-n} (see for instance
\cite[section 6]{c3}).
It is well known that systems of
NLS-NLS time-dependent equations have applications in nonlinear Optics, Hartree-Fock theory for Bose-Einstein
 condensates, among other physical phenomena; see for instance the earlier mathematical works
 \cite{akanbook,ac1,ac2,acr,bt,linwei,mmp,pomp,sirakov}, the more recent list (far from complete)
 \cite{chen-zou,c1,co-oli-tav,liu-wang,soave-tavares} and references therein.
See also a close related work; \cite{c-fract}, in which was studied a close system of coupled NLFS equations.


Here
we are interested in system \eqref{eq:Main}, consisting of coupled NLS equations involving
the so called fractional Laplacian operator (or fractional Schr\"odinger
operator, $(-\D)^s+\l\, \mbox{Id}$).

Note that in dimension $n=1$, \eqref{eq:Main} can also be seen as a system of coupled NLFS-FKdV equations.
In this  case, \eqref{eq:Main} is the corresponding stationary system when one looks for travelling-wave solutions of
the following time-dependent system,
\begin{equation}\label{fractional_NLS-KdV}
\left\{\begin{array}{rcl}
if_t - A_s\,f + |f|^2f+ \b fg & = &0\\
g_t -(A_s\,g)_x +gg_x  + \frac 12\b(|f|^2)_x  & = & 0,
\end{array}\right.
\end{equation}
where $A_s$ stands for the nonlocal fractional Laplacian $(-\Delta)^s$ in dimension $n=1$.

While for $n= 1,\,2,\,3$,  \eqref{eq:Main} can be seen as the stationary system when one looks for standing wave solutions
of the following time-dependent system of coupled NLFS equations,
\begin{equation}\label{fractional_NLS-NLS}
\left\{\begin{array}{rcl}
if_t - (-\D)^s f + |f|^2f+ \b fg & = &0\\
ig_t -(-\D)^s g + \b |f|^2  & = & 0.
\end{array}\right.
\end{equation}

The main goal of this manuscript is to demonstrate that for any $\b>0$, problem
\eqref{eq:Main} has a positive radially symmetric ground state $\wt{\bu}=(\wt{u},\wt{v})\in \Wan\times \Wan$; see Theorems \ref{th:1}, \ref{th:ground2}.

Notice that,  for any $\b\in \R$,
\eqref{eq:Main}  has a unique {\it semi-trivial} positive radially symmetric solution, that we denote by
$\bv_2=(0,V_2)$, where $V_2(x)$ is the unique positive radially symmetric   ground state of
$-(\D)^s v+\l_2 v=\frac 12 v^2$
in $\Wan$; \cite{fl,fls}. Since we are interested in positive ground states,
then we have to show that they are different from the semi-trivial solution $\bv_2$. To do so,
we will demonstrate some properties of the semi-trivial solution which will allow us to show that $\bv_2$ is not a ground state.
For example, we will show that there exists a constant $\L>0$ such that
   for $\b>\L$, $\bv_2$ is a saddle point of the associated energy functional constrained on the corresponding Nehari Manifold,
   which actually is a natural restriction. When
$\b<\L$ then  $\bv_2$ is a strict local minimum of the energy functional on the Nehari Manifold. In this case,
we exclude that $\bv_2$ is a ground state by the construction of a function in the Nehari Manifold with energy
lower than the energy of $\bv_2$.
Precisely, we will demonstrate that there exists a positive radially symmetric ground state of \eqref{eq:Main},
$\wt{\bu}\neq \bv_2$,
either:
$\b>\L$ (see Theorem \ref{th:1}) or $0<\b\le \L$ and $\l_2$ large enough  (see Theorem \ref{th:ground2}).

 The paper is organized as follows. In Section \ref{sec:not-prel} we introduce notation and preliminaries, dealing with
 some background on the fractional Laplacian and we give the definition of ground state.  Section \ref{sec:kr}
 contains some results on the method of the {\it natural constraint} and the main properties about the semi-trivial solution $\bv_2$, that
 we will use in the proof of  the main existence
 results stated and proved in Section \ref{sec:pf2}.
 Finally, in  Section \ref{sec:final} we study the existence of ground states for some systems with an arbitrary number of coupled equations.

\section{Preliminaries and Notation}\label{sec:not-prel}
The nonlocal fractional Laplacian operator $(-\D)^s$ in $\Rn$ is
defined on the Schwartz class of functions $g\in \mathcal{S}$
through the Fourier transform,
\begin{equation}\label{fourier}
    [(-\D)^{\frac{\a}{2}}  g]^{\wedge}\,(\xi)=(2\pi|\xi|)^\alpha \,\widehat{g}(\xi),
\end{equation}
or via the Riesz potential, see for example~\cite{Landkof,Stein}.
Note that $s=1$ corresponds to the standard local
Laplacian operator. See also \cite{Laskin,dpv,fl,fls}, where the fractional Schr\"odinger operator ($(-\D)^s+$Id)
is defined and are analyzed some problems dealing with.

There is another way to define this operator. If
$s=1/2$ the square root of the Laplacian acting on a function
$u$ in the whole space $\Rn$, can be calculated as the
normal derivative on the boundary of its harmonic extension to the
upper half-space $\mathbb{R}^{n+1}_+$, this is so-called Dirichlet
to Neumann operator. Caffarelli-Silvestre; \cite{Caffarelli-Silvestre},
  have shown that this operator can be
realized  in a local way by using one more variable and the so
called $s$-harmonic extension.

More precisely, given  $u$  a regular function defined in $\Rn$
we define its $s$-harmonic extension  to the upper half-space
$\R^{n+1}_+$ by $w=\ext (u)$, as the solution to the
problem \be\label{extension} \left\{
\begin{array}
{rcl@{\qquad}l} -\diver(y^{1-2s}\nabla w)&=&0&\mbox{ in }
\R^{n+1}_+\\ [5pt] w&=&u&\mbox{ on }
\Rn\times\{y=0\}.
\end{array}
\right. \ee
The main relevance of the $s$-harmonic extension
comes from the following identity
 \begin{equation}
\lim\limits_{y\to 0^+}y^{1-2s}\dfrac{\partial w}{\partial
y}(x,y)=-\dfrac{1}{\kappa_s}(-\Delta)^{s}u(x),
\label{normalder}\end{equation} where $\kappa_s$ is a positive
constant.
The above Dirichlet-Neumann procedure
\eqref{extension}-\eqref{normalder} provides a formula for the
fractional Laplacian, equivalent to that obtained from
Fourier Transform by \eqref{fourier}. In that case, the
$s$-harmonic extension  and the fractional Laplacian have
explicit expressions in terms of the Poisson and the Riesz
kernels, respectively,
\begin{equation}
  \label{poisson}
\begin{array}{l}
w(x,y)=P_y^{s}*u(x)=
c_{n,s}\,y^{2s}\displaystyle\int_{\Rn}\dfrac{u(z)}{(|x-z|^{2}+y^{2})^{\frac{n+2s}{2}}}\,dz,
\\ [5mm]
(-\Delta)^{s}u(x)= d_{n,s}\,\displaystyle
P.V.\,\int_{\Rn}\frac{u(x)-u(z)}{|x-z|^{n+2s}}\,dz.
\end{array}
\end{equation}
The natural functional spaces are the homogeneous fractional
Sobolev space  $\dot H^{s}(\Rn)$ and the weighted Sobolev
space $X^{2s}(\R^{n+1}_+)$, that can be defined as the
completion of $\mathcal{C}_{0}^{\infty}(\overline{\R^{n+1}_+})$
and $\mathcal{C}_0^{\infty}(\Rn)$, respectively, under the norms
$$
\begin{array}{l}
\displaystyle\|\phi\|^2_{X^{2s}}=\kappa_s\int_{\R^{n+1}_+}
y^{1-2s}|\nabla
\phi(x,y)|^2\,dx dy,\\[12pt]
\displaystyle\|\psi\|^2_{\dot{H}^{s}}=\int_{\Rn}|2\pi\xi|^{2s}|\widehat{\psi}(\xi)|^2\,d\xi
=\int_{\Rn}|(-\D)^{\frac{s}{2}}\psi(x)|^2\, dx,\\
\end{array}
$$
where $\kappa_s$ is the  constant in \eqref{normalder}.
Notice that, the constants in \eqref{poisson} and $\kappa_s$
satisfy the identity $s\, c_{n,s}\kappa_{s}=d_{n,s}$, and their explicit value
can be seen  in
\cite{brandle-colorado-depablo-sanchez}.
\begin{Remark}
The $s$-harmonic extension operator defined by \eqref{extension} is
an isometry between the spaces $\dot H^{s}(\Rn)$ and
$X^{2s}(\R^{n+1}_+)$, i.e.,
\be\label{isometry} \|\var \|_{\dot
H^{s}}= \| E_s(\var)\|_{X^{2s}}\,,\quad \forall\, \var\in\dot
H^{s}(\R^n).
\ee
Even more, we have  the following inequality for
the trace $\tr (w)=w(\cdot , 0)$, \be\label{trace} \| \tr(w)\|_{\dot
H^{s}}\le\| w\|_{X^{2s}}\, , \quad\forall\, w\in X^{2s}(\R^{n+1}_+),
\ee see \cite{brandle-colorado-depablo-sanchez} for more details.
\end{Remark}
Let us introduce
the following notation:
 \begin{itemize}
\item $E=\Wan$, denotes the fractional Sobolev space, endowed with
scalar product and norm
$$
(u\mid v)_j=\intn \left[ (-\D)^{\frac{s}{2}} u
(-\D)^{\frac{s}{2}}  v + \l_j uv\right] dx,\quad \|u\|_j^2=(u\mid
u)_j,\:\: j=1,2;
$$
\item $\E=E\times E$; the elements in $\E$ will be denoted by $\bu
=(u,v)$; as a norm in $\E$ we will take
$\|\bu\|=\|\bu\|_{\E}^2=\|u\|_1^2+\|v\|_2^2$;
\item $X=X^{2s}(\R^{n+1}_+)$, $\mathbb{X}=X\times X;$
\item  for $\bu\in \E$, the
notation $\bu\geq \bo$, resp. $\bu>\bo$, means that $u,\, v\geq 0$, resp. $u, \, v>0$, for all $j=1,2$.
\end{itemize}
\begin{Remark}
If we define
$$
\frac{\partial w}{\partial \nu^s}=
-\kappa_{s}\lim_{y\to 0^+}y^{1-2s}\,\frac{\partial w}{\partial y}\,,
$$
we can reformulate the main problem \eqref{eq:Main} as
\be\label{eq:sistema-fractional-r2+}\left\{
 \begin{array}{ll}
-\div(y^{1-2s}\nabla w_1) & =0 \quad \mbox{in } \R^{n+1}_+ \\
-\div(y^{1-2s}\nabla w_2) & =0 \quad \mbox{in } \R^{n+1}_+\\
\dfrac{\partial w_1}{\partial \nu^s} +\l_1 w_1 & = w_1^{3}+\b  w_1w_2\quad \mbox{on }\mathbb{R}\times\{y=0\} \\
[2mm]\dfrac{\partial w_2}{\partial \nu^s} +\l_2 w_2 & =
\frac 12 w_2^{2}+\frac 12\b w_1^2\quad \mbox{on }\mathbb{R}\times\{y=0\},
\end{array}\right.
\ee
with $\bw=(w_1,w_2)\in \mathbb{X}$.

Note that if $\bw\in \mathbb{X}$ is solution of \eqref{eq:sistema-fractional-r2+}, then
$\tr (\bw(x,y))=\bw(x,0)\in \E$ is a solution of \eqref{eq:Main}, or equivalently, if $\bu\in \E$ is a solution of \eqref{eq:Main},
then $\ext (\bu)\in \mathbb{X}$
is a solution of \eqref{eq:sistema-fractional-r2+}.

The introduction of this problem is only for the interested reader.
As we will see along the paper, it is not necessary to make use of problem
\eqref{eq:sistema-fractional-r2+}, i.e., all the results for \eqref{eq:Main} are going to be proved without using the
$s$-harmonic extension to the upper half space, $E_s (\cdot)$.
\end{Remark}
For $\bu=(u,\, v)\in \E$,  we set
\be\begin{array}{rcl}\label{eq:II}
I_1(u) &=& \dyle\tfrac 12 \intn (|(-\D)^{\frac{s}{2}}  u|^2+\l_1 u^2)dx -\tfrac{1}{4}\, \intn u^4dx,\\
[4mm]
I_2(v) & = &\dyle\tfrac 12 \intn (|(-\D)^{\frac{s}{2}}  v|^2+\l_2 v^2)dx -\tfrac{1}{6}\, \intn v^3dx,\\
\end{array}
\ee
$$
\Phi (\bu)= I_1(u)+I_2(v)- \tfrac 12\b \intn u^2v\,dx.
$$
We also write
$$
G_\b (\bu)=\tfrac 14\, \intn u^4dx+\tfrac 16\, \intn v^3dx+ \tfrac 12\b \intn u^2v\,dx,
$$
and using this notation we can rewrite the energy functional as
$$
\Phi (\bu)=\frac 12\|\bu \|^2 -G_\b (\bu),\quad \bu\in \E.
$$
We observe that $G_\b$ makes
sense because $\frac n4<s<1\Rightarrow 2^*_s>4$ which implies
the continuous Sobolev embedding $E\hookrightarrow L^{4}(\Rn)$.
Even more, any critical point $\bu\in \E$ of $\Phi$,  gives rise to a solution of
\eqref{eq:Main}.
\begin{Definition}\label{def:ac}
A non-negative critical point  $\wt{\bu}\in \E\setminus\{\bo\}$ is called a ground state of \eqref{eq:Main}
 if its energy $\Phi(\wt{\bu})$ is minimal among all the non-trivial critical points of $\Phi$.
\end{Definition}
\section{The Nehari manifold and properties of $\bv_2$}\label{sec:kr}
Let us set
$$
\Psi(\bu)=(\n\Phi(\bu)|\bu)=(I_1'(u)|u)+(I_2'(v)|v)-\frac 32\,\b \intn u^2v\,dx.$$
We define the  Nehari manifold  by
$$
\cN =\{ \bu\in \E\setminus\{\bo\}: \Psi (\bu)=0\}.
$$
Then, one has that
\begin{equation}\label{eq:gamma}
(\n \Psi(\bu) \mid \bu)= - \|\bu \|^2-\intn u^4\,dx<0\qquad\forall\, \bu\in \cN,
\end{equation}
thus $\cN$ is a smooth manifold locally near any point $\bu\not= \bo$ with $\Psi(\bu)=0$. Moreover,
$\Phi''(\bo)= I_1''(0)+I_2''(0)$ is positive definite, so we infer that $\bo$ is a strict local minimum for $\Phi$. As a consequence,
$\bo$ is an isolated point of the set $\{\Psi(\bu)=0\}$, proving that $\cN$ is a smooth complete manifold of codimension $1$,
and on the other hand there exists a constant $\rho>0$ so that
\be\label{eq:bound1}
\|\bu\|^2>\rho\qquad\forall\,\bu\in \cN.
\ee
Furthermore, by \eqref{eq:gamma} and \eqref{eq:bound1} we can show that $\bu\in \E\setminus\{\bo\}$
is a critical point of $\Phi$ if and only if $\bu\in\cN$ is a critical point of $\Phi$ constrained on $\cN$.

As a consequence, we have the following.
 \begin{Lemma}\label{pr:ac}
 $\bu \in \E$ is a non-trivial critical point of $\Phi$ if and only if
 $\bu\in \cN$ and is a constrained critical point of $\Phi$ on $\cN$.
 \end{Lemma}
\begin{remarks}\label{rem:obs1}
\begin{itemize}
\item[(i)] By the previous arguments, the Nehari manifold $\cN$ is a natural constraint of $\Phi$.
Also, it is relevant to point out that working on the Nehari manifold,
the functional $\Phi$ satisfies the following expression,
\be\label{eq:restriction0}
\Phi|_{\cN}(\bu)= \frac 16\|\bu\|^2+\frac{1}{12}\intn u^4dx=:F(\bu),
\ee
then using  \eqref{eq:bound1} into \eqref{eq:restriction0} we
obtain
\begin{equation}\label{eq:restriction}
\Phi(\bu)\ge   \frac 16\|\bu\|^2>\frac 16 \rho\qquad \forall\, \bu\in\cN.
\end{equation}
Therefore, by \eqref{eq:restriction} the functional $\Phi$ is bounded from below  on $\cN$, as a consequence we will
minimize it on the Nehari manifold. To do so, a remark about compactness is  in order.
\item[(ii)] Analyzing the Palais-Smale (PS) condition, we remember that
working on the radial setting, $H=E_{radial}$, the embedding of $H$ into $L^4(\R^n)$ is compact for $n=2,\, 3$,
but in dimension $n=1$, the embedding of $E$ or $H$ into $L^q(\R)$ for $2<
q<2^*_s$ is not compact; see \cite[Remarque I.1]{Lions-JFA82}.
However, we will analyze all the dimensional cases $n=1,\, 2,\, 3$, proving
that for a PS sequence of $\Phi$ on $\cN$, we can find a subsequence for which the
weak limit is non-trivial and it is a solution of \eqref{eq:Main}. This fact jointly with some properties of
the Schwarz symmetrization will allow us to demonstrate the existence of
positive radially symmetric ground states to \eqref{eq:Main}.
Notice that one could also try to work in the cone of non-negative radially decreasing functions,
where one has the required compactness, in the one-dimensional case, thanks to Berestycki and Lions \cite{BL}, but this is not our approach.
\end{itemize}
\end{remarks}
\begin{Remark}\label{rem:non-trivial}
It is known; \cite{fl,fls},  that the
equation
\be\label{eq:soliton-s} (-\D)^{s} v +v=
v^2,
\ee
with $v\in E$, $v\not\equiv 0$, has a unique  radially symmetric and positive solution,
that we will denote by $V$. Indeed $V$ is a non-degenerate ground state of \eqref{eq:soliton-s} in $H$.

Clearly,  for every $\b\in \R$, \eqref{eq:Main}
already possesses a {\it semi-trivial} solution given by
 $$
 \bv_2 = (0,V_2),
 $$
 where
\be\label{reescale_v2}
V_2(x)=2\l_2 V(\l_2^{1/{2s}}\,x)
\ee
is the unique positive radially symmetric solution of $(-\D)^s v+\l_2
v=\frac 12 v^2$ in $E$.
\end{Remark}
In order to study some useful properties of $\bv_2$, we define
 define the corresponding Nehari manifold associated to $I_2$ in \eqref{eq:II},
$$
\cN_2 =\left\{v\in E : (I_2'(v)|v)=0\right\}=\left\{v\in E : \|v\|_2^2 -\frac 12\intn v^3dx=0\right\}.
$$
Let us denote $T_{\bv_2}\cN$ the tangent space  to  $\cN$ on $\bv_2$. Since
$$
\bh=(h_1,h_2)\in  T_{\bv_2}\cN \Longleftrightarrow
(V_2|h_2)_2= \frac 34\intn V_2^2h_2\,dx ,
$$
it follows that
\begin{equation}\label{eq:tang1}
(h_1,h_2)\in T_{\bv_2} \cN  \Longleftrightarrow h_2\in T_{V_2} \cN_2.
\end{equation}
Then we prove the following.
\begin{Proposition}\label{prop:5}
There exists $\L>0$ such that:
\begin{itemize}
\item[(i)]  if $\b< \L$, then $\bv_2$ is a strict minimum of $\Phi$ constrained on $\cN$,
\item[(ii)] for any   $\b>\L$, then $\bv_2$ is a saddle point of $\Phi$ constrained on $\cN$ with $\dyle\inf_{\cN}\Phi <\Phi(\bv_2)$.
\end{itemize}
\end{Proposition}
\begin{pf}
First, we
observe that if $D^2\Phi_{\cN}$ denotes the second
derivative of $\Phi$ constrained on $\cN$. Using that $\Phi'(\bv_2)=0$ we have that $D^2\Phi_{\cN}
(\bv_2)[\bh]^2=\Phi'' (\bv_2)[\bh]^2$ for all $\bh\in
T_{\bv_2}\cN$.

\

$(i)$ We define
\be\label{eq:Lambda}
\L=\inf_{\varphi\in E\setminus\{
0\}}\frac{\|\varphi\|_1^2}{\intn V_2\varphi^2dx}.
\ee
We have that for $\bh\in  T_{\bv_2}\cN$,
\be\label{eq:Phi-segunda}
\Phi''(\bv_2)[\bh]^2 =\|h_1\|_1^2 +I_2''(V_2)[h_2]^2-\b\intn V_2
h_1^2dx.
\ee
Let us take $\bh=(h_1,h_2)\in T_{\bv_2}\cN$, by
\eqref{eq:tang1} $h_2\in T_{V_2} \cN_2$, then using that $V_2$ is
the minimum of $I_2$ on $\cN_2$, there exists a  constant $c>0$ such
that
\be\label{eq:minimo-pos}
I_2'' (V_2)[h_2]^2\ge c\|h_2\|_2^2.
\ee
Due to \eqref{eq:minimo-pos} jointly with
\eqref{eq:Phi-segunda}, for $\b<\L$, there exists another constant
$c_1>0$ such that, \be\label{minimo} \Phi''(\bv_2)[\bh]^2 \ge
c_1(\|h_1\|_1^2 +\| h_2\|^2), \ee
which proves that  $\bv_2$ is a  strict local minimum of $\Phi$ on $\cN$.

\

$(ii)$
According to \eqref{eq:tang1}, $\bh=(h_1,0)\in T_{\bv_2}\cN$ for any $h_1\in E$.
We have that, for $\b>\L$, there exists $\wt{h}\in E$ with
$$
\L< \frac{\|\wt{h}\|_1^2}{\intn V_2\wt{h}^2dx}<\b,
$$
thus, taking $\bh_0=(\wt{h},0)\in T_{\bv_2}\cN$, by \eqref{eq:Phi-segunda} we find
\be\label{eq:*1}
\Phi''(\bv_2)[\bh_0]^2 =\|\wt{h}\|_1^2 -\b\intn V_2 \wt{h}^2dx<0.
\ee
On the other hand,
by \eqref{eq:tang1}, and using again that $V_2$ is
the minimum of $I_2$ on $\cN_2$, we have that there exists $c>0$ such that
$$
I_2''(V_2)[h]^2\ge c\|h\|_2^2, \:\forall\, h\in T_{V_2}\cN_2.
$$
Finally, by \eqref{eq:Phi-segunda}, $\Phi''(\bv_2)[(0,h)]^2=I_2''(V_2)[h]^2$ for any $h\in T_{V_2}\cN_2$.
Thus we have that there exists a  constant $c>0$ such that
$$
\Phi''(\bv_2)[\bh_1]^2\ge c\|\bh_1\|^2,\quad \forall\,\bh_1=(0,h_1)\in T_{\bv_2}\cN.
$$
\end{pf}
\section{Ground state solutions}\label{sec:pf2}
The first result on the existence of ground states is given for the coupling parameter $\b>\L$ in the following.

\begin{Theorem}\label{th:1}
Assume  $\b>\Lambda$, then $\Phi$  has a positive radially symmetric ground state $\wt{\bu}$, and there holds
$\Phi(\wt{\bu})<\Phi(\bv_2)$.
\end{Theorem}
\begin{pf}
By the Ekeland's variational principle;
\cite{eke}, there exists a PS sequence $\{\bu_k\}_{k\in
\mathbb{N}}\subset\cN$, i.e.,
\be\label{eq:PS1}
\Phi (\bu_k)\to
c_\cN=\inf_{\cN}\Phi
\ee
\be
\n_{\cN}\Phi(\bu_k)\to 0.
\ee
By \eqref{eq:restriction0} and \eqref{eq:PS1}, we find that $\{\bu_k\}$ is a
bounded sequence on $\E$, hence for a subsequence, we can assume that
\be\label{eq:weakly}
\bu_k\rightharpoonup \bu_0\quad\mbox{weakly in } \E,
\ee
\be\label{eq:local}
\bu_k\to \bu_0\quad \mbox{strongly
in }\mathbb{L}^q_{loc}(\R)=L^q_{loc}(\R)\times  L^q_{loc}(\R)\quad \forall\, 1\le q<2^*_s,
\ee
and also $\bu_k\to \bu_0$ a. e. in $\Rn$. Since $\cN$ is closed we have that $\bu_0\in\cN$, even more, using that $\bo$ is an isolated point
the set $\{ \Psi (\bu)=0\}$ we infer that $\bu_0\neq \bo$.
On the other hand,
the constrained gradient satisfies
\be\label{eq:to0}
\n_{\cN}\Phi (\bu_k)=\Phi' (\bu_k)-\eta_k \Psi'(\bu_k)\to 0,
\ee
where $\eta_k$ is the
corresponding Lagrange multiplier. Taking the scalar product with
$\bu_k$ in \eqref{eq:to0}, since $\bu_k\in\cN$ we have that  $(\Phi'(\bu_k)\mid
\bu_k)=\Psi(\bu_k)=0$, then we infer that $\eta_k (\Psi'(\bu_k)\mid
\bu_k)\to 0$; this jointly with
\eqref{eq:gamma},\eqref{eq:restriction0} and the fact that $\|\Psi'(\bu_k)\|\leq C<\infty$ imply  that $\eta_k\to 0$ and therefore
$\Phi'(\bu_k)\to 0$.

As a consequence of the discussion above, although we do not know that $\bu_k\to\bu_0$ in $\E$, we infer that $\bu_0\in\cN$
is a non-trivial critical point of $\Phi$ and by Lemma \ref{pr:ac} it is also a non-trivial critical point of $\Phi$ on $\cN$.

Moreover, using that  $\bu_0\in\cN$ jointly with \eqref{eq:restriction0} and the Fatou's Lemma,  we find
$$
\begin{array}{rcl}
\Phi (\bu_0) & = & \dyle F(\bu_0)\\
& \le & \dyle\liminf_{k\to\infty} F(\bu_k)\\
 & = & \dyle\liminf_{k\to\infty}\Phi(\bu_k)= c_\cN.
 \end{array}
$$

As a consequence, $\bu_0$ is a least energy solution of \eqref{eq:Main}.
By Proposition \ref{prop:5}-$(ii)$ we know that
necessarily $\Phi(\bu_0)<\Phi(\bv_2)$. Additionally, by the maximum principle in the fractional setting; \cite{cabre-sire}, applied to
the second equation in \eqref{eq:Main}, we have that
$v_0>0$.
In order to show that also  $u_0>0$,
first we prove the following.

\

{\it Claim.} We can assume without loss of generality that ${u}_0\ge 0$.

\

To prove it, we consider
$|\bu_0|=(|{u}_0|,v_0)$, then we have two cases:
\begin{enumerate}
\item If $|\bu_0|\in\cN$, by the Stroock-Varopoulos inequality;
\cite{stroock, varop},
\be\label{eq:stroock-var}
\|(-\Delta)^\frac{s}{2} (|u|)\|_{L^2}\le
\|(-\Delta)^\frac{s}{2} (u)\|_{L^2},
\ee
we have, in particular, that
 $\||u|\|_1\le \| u\|_1$, then we obtain
$$
\Phi (|\bu_0|)\le \Phi(\bu_0)=c_\cN.
$$
Then, by similar arguments as in \cite[Theorem 4.3]{Willem},
we have that $|\bu_0|$ is a non-negative ground state.
\item If $|\bu_0|\not\in\cN$, we take the unique $t>0$, $t\neq 1$ such that
$t|\bu_0|\in\cN$, which comes from
\be\label{eq:tpos}
\|\, |\bu_0|\, \|^2=t^2\intn u_0^4dx+t\left(
\frac 12\intn v_0^3\,dx+\frac 32\b \intn
{u}_0^2{v}_0 \,dx\right).
\ee
Since $\bu_0\in\cN$, then
\be\label{eq:tpos2}
\|\, \bu_0\, \|^2=\intn {u}_0^4dx+
\frac 12\intn {v}_0^3\,dx+\frac 32\b \intn
{u}_0^2{v}_0 \,dx.
\ee
By \eqref{eq:tpos}, \eqref{eq:tpos2} and again the Stroock-Varopoulos inequality \eqref{eq:stroock-var}, we infer that
\be\label{eq:tpos3}
\begin{array}{rcl}
 & & \dyle t^2\intn {u}_0^4dx+t\left(
\frac 12\intn {v}_0^3\,dx+\frac 32\b \intn
{u}_0^2{v_0} \,dx\right)\\
[4mm]
& & \le
\dyle\intn {u}_0^4dx+
\frac 12\intn {v}_0^3\,dx+\frac 32\b \intn
{u}_0^2{v}_0 \,dx.
\end{array}
\ee
Using that $t\neq 1$, as a consequence of \eqref{eq:tpos3} we deduce that $0<t<1$ and the inequality in \eqref{eq:tpos3} is strict.
Hence, by \eqref{eq:restriction0} jointly with \eqref{eq:stroock-var} and $t<1$  we obtain
$$
\begin{array}{rcl}
\Phi (t|\bu_0|) & = &\dyle t^2\|\, |\bu_0|\, \|^2 +t^4\frac{1}{12}\intn u_0^4dx\\
[3mm]
& < & \dyle \|\, |\bu_0|\, \|^2 +\frac{1}{12}\intn u_0^4dx\\
[3mm]
& \le & \Phi(\bu_0)=c_\cN.
\end{array}
$$
This is a contradiction because  $t|\bu_0|\in\cN$. Therefore $|\bu_0|\in\cN$ and the claim is proved.
\end{enumerate}
Once we can assume without loss of generality that ${u}_0\ge 0$, by the maximum principle applied to the first equation in \eqref{eq:Main}
we find ${u}_0> 0$ proving that indeed $\bu_0$ is a positive ground state.

To finish the proof, we have to show that the ground state is indeed radially symmetric.

If $\bu_0$ is not radially symmetric, we set $\wt{\bu}=\bu_0^\star=(u_0^\star,v_0^\star)$,
where $u_0^\star,\,{v}_0^\star$ denote the Schwarz symmetric functions associated to
$u_0,\, v_0$ respectively. By the properties of the
Schwarz symmetrization; see for instance \cite{fmm} for the fractional setting and
\cite{ban} for the classical one, there hold
\be\label{eq:primera}
\|{\bu}^\star\|^2\le \|{\bu}\|^2, \qquad
G_\b({\bu}^\star)\ge G_\b({\bu}).
\ee
Furthermore,  there exists a unique $t_\star>0$ such that $t_\star\,\wt{\bu}\in
{\mathcal{N}}$. If $t_\star=1$, by \eqref{eq:primera} we have $\Phi (\wt{\bu})\le \Phi (\bu_0)=c_\cN$ with $\wt{\bu}\in\cN$
thus $\wt{\bu}$ is a positive radially symmetric ground state of \eqref{eq:Main}.

On the contrary, i.e., if $t_\star\neq 1$, as in \eqref{eq:tpos}, $t_\star$ comes from
\be\label{eq:t}
\| \wt{\bu}\|^2=
t_\star^2\intn (u_0^\star)^4dx+t_\star\left( \frac 12\intn
(v_0^\star)^3dx+\frac 32\b \intn (u_0^\star)^2v_0^\star
\,dx\right).
\ee
Due to $\bu_0\in\cN$,
\eqref{eq:primera}, \eqref{eq:t}, the fact that $\bu_0>\bo$ and
$t_\star>0$ we find
\be\label{t_star}
\begin{array}{rcl}
& & \dyle \intn u_0^4\,dx+ \frac 12\intn v_0^3\,dx+\frac 32\b \intn u_0^2v_0\,dx \\ & & \\
& \ge &\dyle  t_\star^2\intn u_0^4dx+t_\star\left( \frac 12\intn
v_0^3\,dx+\frac 32\b \intn u_0^2v_0 \,dx\right).
\end{array}
\ee
Thus, using that $0<t_\star\neq 1$ in \eqref{t_star}, we obtain  $0<t_\star< 1$, this and \eqref{eq:primera} show that
\be\label{eq:previa}
\Phi(t_\star\,\wt{\bu})= \frac 16
t_\star^2\|{\bu}^\star\|^2+\frac{1}{12}t_\star^4\intn (u_0^\star)^4\,dx <
\frac 16 \|\bu_0\|^2+\frac{1}{12}\intn u_0^4\,dx=\Phi(\bu_0)=c_\cN,
\ee
with $t_\star\wt{\bu}\in\cN$ which is a contradiction with \eqref{eq:previa}, proving that $t_\star=1$ and as above, then we finish the proof.
\end{pf}
The second result about existence of ground states cover the range $0<\b\le\L$, provided  $\l_2$ is large enough.
\begin{Theorem}\label{th:ground2}
There exists $\L_2>0$ such that if $\l_2>\L_2$, System \eqref{eq:Main} has a radially symmetric ground state
$\wt\bu>\bo$ for every $0<\b\le\L$.
\end{Theorem}
\begin{pf} Arguing in the same way as in the proof of Theorem \ref{th:1},
we prove that there exists a radially symmetric ground state
$\wt{\bu}\ge \bo$. Moreover, in Theorem \ref{th:1} for $\b>\L$ we
proved that $\wt{\bu}>\bo$. Now we need to show that for $0<\b\le \L$
indeed $\wt{\bu}>\bo$ which follows by the maximum principle
provided $\wt{\bu}\neq \bv_2$. Taking into account Proposition
\ref{prop:5}-$(i)$, $\bv_2$ is a strict local minimum of $\Phi$ on $\cN$,  and this does not guarantee
that $\bu_0\not \equiv\bv_2$. Following \cite{c3}, the  idea
consists on the construction of a function
$\bu_0=(u_0,v_0)\in\cN$ with $\Phi(\bu_0)<\Phi(\bv_2)$. To do so,
since $\bv_2=(0,V_2)$ is a local minimum of $\Phi $ on $\cN$
provided $0<\b<\L$, we cannot find $\bu_0$ in a neighborhood of
$\bv_2$ on $\cN$. Thus, we define $\bu_0=t(V_2,V_2)$ where  $t>0$ is
the unique value such that $\bu_0\in \cN$.

Now, we will show that
$$
\bu_0=t(V_2,V_2)\in\cN\quad \hbox{with}\quad \Phi(\bu_0)<\Phi(\bv_2),
$$
for $\lambda_2$ large enough.

Notice that $t>0$ comes from $\Psi (\bu_0)=0$, i.e.,
\be\label{eq:condicion1}
t^2\|(V_2,V_2)\|^2-t^4\intN V_2^4\,dx
-\tfrac{1}{2}t^3(1+3\beta)\intN V_2^3\,dx=0.
\ee
We also have
\be\label{eq:norma doble}
\|(V_2,V_2)\|^2=2\|V_2\|_2^2+(\lambda_1-\lambda_2)\intN V_2^2\,dx.
\ee Moreover, since $V_2\in\cN_2$, we have \be\label{eq:norma
simple} \|V_2\|_2^2-\tfrac{1}{2}\intN V_2^3\,dx=0.
\ee
Substituting \eqref{eq:norma doble}  and \eqref{eq:norma simple} in
\eqref{eq:condicion1} it follows
\be\label{eq:condicion2}
t^2\left(\intN V_2^3\,dx+(\lambda_1-\lambda_2) \intN
V_2^2\,dx\right) -t^4\intN V_2^4\,dx -\tfrac{1}{2}t^3(1+3\beta)\intN
V_2^3\,dx=0.
\ee
Hence, applying the scaling \eqref{reescale_v2}  yields
\be\label{cambio p}
\intN V_2^r\,dx=2^r\lambda_2^{r-\frac{n}{2s}}\intN V^r\,dx.
\ee
 Subsequently, substituting \eqref{cambio p} for $r=2,3,4$ into
\eqref{eq:condicion2} and dividing by $2^3\lambda_2^{3-\frac{n}{2s}}t^2$
we have that
\be\label{eq:condicion}
\intN V^3\,dx+\dfrac{\lambda_1-\lambda_2}{2\lambda_2}\intN
V^2\,dx-2\lambda_2t^2\intN V^4\,dx-\tfrac{1}{2}t(1+3\beta)\intN
V^3\,dx=0.
\ee
Moreover, by \eqref{eq:restriction0},
\eqref{eq:norma doble} and \eqref{eq:norma simple} we find
respectively the expressions
\be\label{eq: forma de Phi(w)}
\Phi(\bu_0)=\tfrac{1}{6}t^2\left( \intN
V_2^3\,dx+(\lambda_1-\lambda_2) \intN V_2^2\,dx\right)
+\tfrac{1}{12}t^4\intN V_2^4\,dx,
\ee
\be\label{eq: forma de Phi(bv_2)}
\Phi(\bv_2)=I_2(V_2)=\tfrac{1}{2}\|V_2\|_2^2-\tfrac{1}{6}\intN
V_2^3=\tfrac{1}{12}\intN V_2^3.
\ee
By \eqref{eq: forma de Phi(w)}, \eqref{eq: forma de Phi(bv_2)} we have $\Phi(\bu_0)< \Phi(\bv_2)$ is equivalent to
\be\label{40}
\tfrac{1}{6}t^2\left( \intN
V_2^3\,dx+(\lambda_1-\lambda_2) \intN V_2^2\,dx\right)
+\tfrac{1}{12}t^4\intN V_2^4\,dx- \tfrac{1}{12}\intN V_2^3\,dx <0,
\ee
and then, applying again \eqref{cambio p} and multiplying
\eqref{40} by $6\lambda_2^{\frac{n}{2s}-3}$, we actually have
\be\label{desigualdad de w}
t^2\left( \intN V^3\,dx+\dfrac{\lambda_1-\lambda_2}{\lambda_2}\intN V^2\,dx\right)
+\tfrac{1}{2}t^4\lambda_2\intN V^4\,dx-\tfrac{1}{2}\intN V^3\,dx<0.
\ee
For $\l_2$ large enough we find that \eqref{eq:condicion} will provide us with \eqref{desigualdad de w}.
Therefore, there exists a positive constant $\L_2$ such that for $\l_2>\L_2$ inequality \eqref{desigualdad de w} holds,
proving that
$$
\Phi(\wt{\bu})\leq\Phi(\bu_0)< \Phi(\bv_2).
$$
Finally, this  shows that $\wt{\bu}\neq \bv_2$ and we finish.
\end{pf}
\section{Systems with more than $2$ equations}\label{sec:final}
In this last subsection, we deal with some extended systems of \eqref{eq:Main} to more than
two equations.

We start with the study of the
following system coming from NLFS-2FKdV equations if $n=1$ or 3NLFS equations if $n=1,\, 2,\, 3$,
\begin{equation}\label{NLS-KdV2-3}
\left\{\begin{array}{rcl}
(-\D)^s u +\l_0 u & = & u^3+\b_1 uv_1+\b_{2}uv_2,\\
(-\D)^s v_1 +\l_1 v_1 & = & \frac 12 v_1^2+\frac 12\beta_{1} u^2,\\
(-\D)^s v_2+\l_2 v_2 & = & \frac 12v_2^2+\frac 12\beta_{2} u^2,
\end{array}\right.
\end{equation}
where $u,\, v_1,\,v_2\in E$.
This system can be seen as a perturbation of \eqref{eq:Main} when  $|\b_{1}|$ or $ |\b_{2}|$ is small.

We use similar notation as in previous sections  with natural
meaning, for example, $\E= E\times E \times
E$, $\bo=(0,0,0)$, \be\label{eq:Phi3} \Phi (\bu)=\frac 12 \|
\bu\|^2-\frac 14\intn u^4\, dx-\frac 16\intn (v_1^3+v_2^3)\,dx-\frac
12\intn u^2(\b_1v_1+\b_2v_2)\,dx \ee \be\label{eq:N3} \cN=\{\bu\in
\E\setminus\{\bo\}:  (\Phi'(\bu)| \bu)=0\}, \ee etc.

Let $U^*,\, V_j^*$ be the unique positive radially symmetric solutions of $(-\D)^s u+\l_0 u=u^3$, $(-\D)^s v+\l_j v=\frac 12 v^2$
in $E$ respectively, $j=1,\,2$; see \cite{fl,fls}.
\begin{Remark}\label{rem:semi-trivial3}
The unique non-negative semi-trivial solutions of \eqref{NLS-KdV2-3} are given by
$\bv_1^*=(0,V_1^*,0)$, $\bv_2^*=(0,0,V^*_2)$ and $\bv_{12}^*=(0,V_1^*,V^*_2)$.
\end{Remark}
Following Section \ref{sec:pf2}, the first result about existence of ground states is the following.
\begin{Theorem}\label{th:3}
Assume $\b_{j}>\L_j$ for $j=1,\, 2$, then \eqref{NLS-KdV2-3} has a
positive radially symmetric ground state $\wt{\bu}$.
\end{Theorem}
\begin{pf}
We define
\be\label{eq:Lambda-j}
\L_j=\inf_{\varphi\in E\setminus\{
0\}}\frac{\|\varphi\|_0^2}{\intn V_{j}^*\varphi^2dx}\qquad j=1,2.
\ee
where $\|\cdot\|_0$ is the norm in $E$ with $\l_0$.

As in Proposition \ref{prop:5}-$(ii)$, using that $\b_{j}>\L_j$,
j={1,\,2}, one can show that both $\bv_1^*$, $\bv_2^*$ are saddle
points of the energy functional $\Phi$  (defined by \eqref{eq:Phi3})
constrained on the Nehari manifold $\cN$ (defined by
\eqref{eq:Phi3}). Then
\be\label{eq:energy3}
c_\cN=\inf_{\cN}\Phi<\min\{ \Phi(\bv_1^*),\Phi
(\bv_2^*)\}<\Phi(\bv_{12}^*)=\Phi(\bv_1^*)+\Phi (\bv_2^*).
\ee
By the Ekeland's variational principle, there exists a PS sequence
$\{\bu_k\}_{k\in \mathbb{N}}\subset\cN$, i.e.,
\be\label{eq:PS1-3}
\Phi (\bu_k)\to c_\cN \ee \be \n_{\cN}\Phi(\bu_k)\to 0.
\ee
The lack of compactness can be circumvent arguing in a similar way as in the proof of Theorem \ref{th:1},
proving that for a subsequence,
$\bu_k\rightharpoonup \wt{\bu}$ weakly in $\E$ with $\wt{\bu}\gvertneqq\bo$, $\wt{\bu}\in\cN$ a critical point of $\Phi$ satisfying
$ \Phi(\wt{\bu})=c_\cN$, then $\wt{\bu}$ is a non-negative ground state.

To prove the positivity of $\wt{\bu}$, if one supposes that the
first component $u^* \equiv 0$, since the only non-negative
solutions of \eqref{NLS-KdV2-3} are the semi-trivial solutions
defined  in Remark \ref{rem:semi-trivial3}, we obtain a contradiction
with \eqref{eq:energy3}. Furthermore, if the second or third
component vanish then $\wt{\bu}$ must be $\bo$, and this is not
possible because $\Phi|_{\cN}$ is bounded bellow by a positive
constant like in \eqref{eq:restriction}, then  $\bo$ is an isolated
point of the set $\{\bu\in \E\, :\: \Psi(\bu)=(\Phi'(\bu)|\bu)=0\}$, proving
that $\cN$ is a complete manifold, as in the previous sections. Then, the maximum principle
shows that $\wt{\bu}>\bo.$ Finally, to show that we have a radially symmetric ground state, we argue as in the proof of Theorem \ref{th:1}.
\end{pf}

Furthermore, following the ideas in the proof of Theorem \ref{th:ground2} we have the following.
\begin{Theorem}\label{th:4} Assume that $\b_{1},\b_{2}>0$ (but not necessarily $\b_j>\L_j$ as in Theorem \ref{th:3}).
Then there exists a positive radially symmetric ground state $\wt\bu$ provided $\l_1,\l_2$ are sufficiently large.
\end{Theorem}
\begin{pf}
The proof follows the same ideas as the one of Theorem
\ref{th:ground2} with appropriate changes. For example, in order to prove the positivity, one has to show that there exists $\bu_0\in\cN$
with $\Phi(\bu_0)<\min\{ \Phi(\bv_1^*),\Phi
(\bv_2^*)\}$, that holds true provided $\l_1,\l_2$ are large enough. We  omit details for
short.
\end{pf}

Plainly we can extend these results to systems with  an arbitrary number of equations $N>3$ as the following,
\be\label{eq:system-N}
\left\{\begin{array}{rcl}
(-\D)^s u +\l_0 u & = & \dyle u^3+\sum_{k=1}^{N-1}\b_{k}\, uv_k \\ & & \\
(-\D)^s v_j +\l_j v_j & = & \dyle\frac 12 v_j^2+\frac 12\b_{j} u^2; \qquad j=1,\cdots,N-1.\\
\end{array}\right.
\ee
Arguing as in  Theorems \ref{th:3}, \ref{th:4} we can show the next result.
\begin{Theorem}\label{th:5}
There exists a positive radially symmetric ground state of \eqref{eq:system-N} if
\begin{itemize}
\item either
$$
\b_{k}>\L_k=\inf_{\varphi\in E\setminus\{
0\}}\frac{\|\varphi\|_0^2}{\intn V_{k}^*\varphi^2dx};\quad
k=1,\cdots N-1,
$$
where $V_k^*$ denotes the unique positive radial solution of $\D v+\l_k v=\frac 12 v^2$ in $E$; $k=1,\cdots, N-1$,
\item or $\b_j>0$ are arbitrary and  $\l_j$ are large enough; $j=1,\ldots,N-1$.
\end{itemize}
\end{Theorem}
\begin{Remark}
As was commented in \cite{c3} for the local setting, here in the nonlocal fractional framework,
another natural extension of \eqref{NLS-KdV2} to more than two
equations  different from \eqref{NLS-KdV2-3} is the following system
coming from 2NLFS-FKdV equations if $n=1$ or 3NLFS equations if $n=1,\, 2,\, 3$, \be\label{eq:NLS2-KdV}
\left\{\begin{array}{rcl}
(-\D)^s u_1 +\l_1 u_1 & = & u_1^3+\b_{12} u_1u_2^2+\b_{13}u_1v \\
(-\D)^s u_2 +\l_2 u_2 & = &  u_2^3+\frac 12\beta_{12} u_1^2u_2+\b_{23} u_2v\\
(-\D)^s v+\l v & = & \frac 12v^2+\frac 12\beta_{13} u_1^2+\frac
12\b_{23}u_2^2.
\end{array}\right.
\ee
We denote $U_j$  the unique positive radially symmetric solution of $(-\D)^s u +\l_j u  = u^3$ in $E$; $j=1,\, 2$;
and $V$ the corresponding positive radially symmetric solution to
 $(-\D)^s v +\l v = \frac 12 v^2$ in $E$.

Note that the non-negative radially symmetric semi-trivial solution $(0,0,V)$ is
a strict local minimum of the associated energy functional
constrained on the corresponding Nehari manifold provided
$$
\b_{j3}<\L_{j}=\inf_{\varphi\in E\setminus\{
0\}}\frac{\|\varphi\|_{\l_j}^2}{\intn V\varphi^2dx}\qquad j=1,2.
$$
While if either $\b_{13}>\L_1$ or $\b_{23}>\L_{2}$ then
 $(0,0,V)$ is a saddle  point of $\Phi $ on $\cN$.

There also exist semi-trivial solutions coming from the solutions
studied in Section \ref{sec:pf2}, with the first component or the
 second one $\equiv 0$. This fact makes different the analysis of
 \eqref{eq:NLS2-KdV} with respect to the previous studied systems \eqref{NLS-KdV2-3} and \eqref{eq:system-N}.
To finish, one could study more general extended systems of
\eqref{NLS-KdV2-3}, \eqref{eq:NLS2-KdV}
 with $N=m+\ell$; coming from $m$-NLFS and $\ell$-FKdV coupled equations with $m,\,\ell\ge 2$ in the one dimensional case, or
 $N$-NLFS equations if $n=1,\,2,\,3$. Indeed, the existence of positive ground states it is still unknown in the local setting ($s=1$)
 for this last kind of systems, including  \eqref{eq:NLS2-KdV} with $s=1$.
\end{Remark}

\enddocument